\documentclass[a4paper,10pt]{article}

\usepackage[latin1]{inputenc}
\usepackage[T1]{fontenc}

\usepackage{psfrag}
\usepackage{amssymb,amsmath,amsfonts,amsthm}

\usepackage{bm, bbm}
\usepackage{graphicx}
\usepackage{epsfig}
\usepackage{natbib}
\usepackage{caption}

\usepackage{hyperref} \hypersetup{dvips, bookmarks, bookmarksnumbered,
  colorlinks=false, breaklinks, pdfborder=0 0 0, pdftitle=Estimation of
  the volume of an excursion set of a Gaussian process using intrinsic
  Kriging, pdfauthor=Emmanuel Vazquez, pdfkeywords=Excursion set
  Gaussian process Intrinsic Kriging Quantile estimation Failure
  probability Design of experiments}

\bibpunct{(}{)}{~;}{a}{,}{,}

\newcommand{\RR}{{\mathbb{R}}}
\newcommand{\NNN}{{\mathbb{N}}}

\newcommand{\NN}{{\mathcal{N}}}
\newcommand{\XX}{{\mathbb{X}}}
\newcommand{\PR}{\mathsf{P}}
\DeclareMathOperator{\var}{var}
\DeclareMathOperator{\cov}{cov}
\DeclareMathOperator{\EE}{E}
\newcommand{\argmin}{\mathop{\rm argmin\;}}

\newcommand{\one}{{\mathbf{1}}}
\newcommand{\tr}{^{\mathsf{T}}}

\newcommand{\x}{{x}}
\newcommand{\y}{{{y}}}

\providecommand{\abs}[1]{\lvert#1\rvert}

\providecommand{\ns}[1]{\lVert#1\rVert}

\renewcommand{\hat}{}
\newcommand{\bmid}{\;\big\lvert\;}

\newtheorem{defn}{Definition}
\theoremstyle{plain}
\newtheorem{prop}{Proposition}

\newtheorem{thm}{Theorem}

\newenvironment{myabstract}
{
\noindent\ignorespaces%
\rule{\textwidth}{.2pt}%
\par%
\setlength{\parindent}{0pt}%
\setlength{\parskip}{1.2ex plus 0.5ex minus 0.2ex}%
{\bf Abstract} \hspace{0.5em} --- \hspace{0.5em}%
}
{
\par
\setlength{\parskip}{0ex plus 0.5ex minus 0.2ex}%
\noindent\ignorespacesafterend%
\rule{\textwidth}{.2pt}%
}

\title{Estimation of the volume of an excursion set of a Gaussian
    process using intrinsic Kriging}

\author{\parbox{0.95\textwidth}{Emmanuel Vazquez$^{\rm a}$ \\
{ \small E-mail address: emmanuel.vazquez@supelec.fr} \\
Miguel Piera Mart\'inez$^{\rm a,b}$ \\
{ \sl \small a: Département Signaux et Systèmes Électroniques, 
      Supélec, 91192  Gif-sur-Yvette, France} \\
{\sl \small b: Laboratoire des Signaux et Systèmes,
      CNRS, Supélec, Université Paris-Sud,  91192  Gif-sur-Yvette, France}
}}
\date{\today}

\begin{document}

\maketitle

  \begin{myabstract}
    Assume that a Gaussian process $\xi$ is predicted from $n$ pointwise
    observations by intrinsic Kriging and that the volume of the
    excursion set of $\xi$ above a given threshold $u$ is approximated
    by the volume of the predictor.  The first part of this paper gives
    a bound on the convergence rate of the approximated volume.  The
    second part describes an algorithm that constructs a sequence of
    points to yield a fast convergence of the approximation.  The
    estimation of the volume of an excursion set is a highly relevant
    problem for the industrial world since it corresponds to the
    estimation of the failure probability of a system that is known only
    through sampled observations.


  \end{myabstract}

%

\section{Introduction}
\vspace{-1ex}

The problem to be considered in this paper is the estimation of the
probability
\begin{equation}
  \mathcal{P}_u := \PR\{ f(X) \geq u\},
  \label{eq:1}
\end{equation}
where $f(\x)$ is a real function defined over an arbitrary set~$\XX$
($\XX=[0,1]^d$ or $\XX=\RR^d$, in most situations) endowed with a
probability measure $\mu$ and $X\in\XX$ is a random vector with
the distribution $\mu$. In practice, the estimation of~(\ref{eq:1}) is
based on a finite sequence of evaluations of $f$ at points
$(x_{i})_{1\leq i \leq n}$ in~$\XX$.
Another way of looking at~(\ref{eq:1}) is via the excursion set
\begin{equation}
	A_u(f) := \{\x \in \XX : f(\x) \geq u\}
	\label{eq:2}
\end{equation}
of the function $f$ above the level $u$, since $\PR\{ f(X) \geq
u\}$ is the volume  $\mu(A_u(f))$, hereafter denoted
by $\abs{A_u(f)}$.

Such a problem is frequently encountered in engineering: the probability
that the inputs of the system will generate a level of a function of the
outputs that exceeds a specified reference level may be expressed
as~(\ref{eq:1}) (where in this case, $X$ is the vector of the inputs
of the system and $f$ is a statistic of the outputs).  Since to obtain
the value of $f$ at a given $\x$ may be very expensive in practice,
because it may involve heavy computer codes for instance, it is often
essential to estimate $\mathcal{P}_u$ using as few evaluations of $f$ as
possible.

To overcome the problem of evaluating $f$ many times, one possible
approach is to estimate $\abs{A_u(f_n)}$ instead of $\abs{A_u(f)}$,
where $\hat f_n$ is an approximation of $f$ constructed from a small set
$\{f(\x_1),\ldots, f(\x_n)\}$ of pointwise evaluations. Such an
approximation can be obtained by assuming that $f$ is a sample path of a
Gaussian random process $\xi$ and by using a linear predictor $\hat
\xi_n$ of $\xi$ constructed from $\xi(x_i)$, $i =1,\ldots,n$. 
In this paper, intrinsic Kriging \citep{Mat73} will be used to obtain
$\hat \xi_n$.  We shall show in Section~\ref{sec:excursion-set-volume}
that this method is likely to give faster convergences than the classical
Monte Carlo estimators, depending on the regularity of $\xi$.

A second step is to choose a sequence of evaluation points $(\x_i)$
so that $\abs{A_u(\xi_n)}-\abs{A_u(\xi)}$ conditioned on the random
variables $\xi(x_i)$, $i\leq n$, converges rapidly to zero.
Section~\ref{sec:SUR} presents an acceleration algorithm based on
computing an upper bound of the mean square error of volume
approximation conditioned on the events $\{\xi(x_i)=f(x_i)$,
$i=1\ldots,n\}$: a point $x_{n+1}$ is selected so that evaluating
$f(\x_{n+1})$ yields the potential largest decrease of the upper bound.
Section~\ref{sec:example}  provides a numerical example.

\section{Excursion set volume estimation by intrinsic Kriging}
\label{sec:excursion-set-volume}

This section deals with the estimation of the
probability~$\mathcal{P}_u$ from observations of $f$ at a finite
sequence of points $(x_i)_{1 \leq i\leq n}$.  As mentioned above,
$\mathcal{P}_u$ is the volume of $A_u(f)$ under the probability
distribution $\mu$. We assume moreover that $f$ is a sample path of a
(separable) Gaussian process $\xi$, with mean $m(x)$, $\x\in\XX$, and
covariance $k(x,y)$, $(x,y)\in\XX^2$.

\subsection{Monte Carlo estimation}

Monte Carlo is a commonly used method to estimate $\abs{A_u(\xi)}$. The
volume of excursion of a Gaussian process $\xi$ may be estimated by
\begin{equation}
\abs{A_{u}(\xi)}_l := \frac{1}{l}\sum_{i=1}^l \one_{\{\xi(X_i)\geq u \}}
\quad {\rightarrow_l} \quad \abs{A_{u}(\xi)} \qquad\mbox{a.s.}
\label{eq:3}
\end{equation}
where the $X_i$s are independent random variables with distribution $\mu$.
The estimator~(\ref{eq:3}) is unbiased, since $E[\abs{A_u(\xi)}_l \mid
\xi]=\abs{A_u(\xi)}$, and 
$$
\EE\big[(\abs{A_{u}(\xi)}_l - \abs{A_{u}(\xi)})^2\mid \xi \big] =
\frac{1}{l}\abs{A_{u}(\xi)}\big( 1- \abs{A_{u}(\xi)}\big)\,.
$$
If evaluating $f$ (a sample path of $\xi$) at many points of $\XX$ is
not particularly demanding, then estimating $\abs{A_u(f)}$ is
straightforward. However, if $\abs{A_{u}(f)}$ is small, then the
variance of the Monte Carlo estimator is approximately
${\abs{A_{u}(f)}}/{l}$. To achieve a given standard deviation $\kappa
\abs{A_{u}(f)}$, with $\kappa>0$ small, the required number of
evaluations is approximately $1/({\kappa^2 \abs{A_{u}(f)}})$, i.e. it is
high.  Thus, the convergence of~(\ref{eq:3}) may be too slow in many
real applications where doing a lot of evaluations of $f$ may not be
affordable (for instance, $f$ may be a complex computer simulation and
may take hours or days to run). Of course, many other methods have been
proposed to improve the basic Monte Carlo convergence. For instance,
methods based on importance sampling, on cross-entropy
\citep{rubinstein99}, on the classical extreme value theory
\citep[e.g.][]{embrechts97:_model_extrem_event_insur_finan}, etc.  They
are not considered here for the sake of brevity.

\subsection{Estimation based on an approximation}

An alternative approach is to replace $f$ by an approximation $\hat f_n$
constructed from a set of $n$ point evaluations of $f$. Provided $\hat
f_n$ converges rapidly enough to $f$, one expects a good estimation of
the excursion sets and their volume using only a few evaluations of $f$.
There are many ways of constructing such an approximation.  Let us
mention two classical methods: regularized regressions in reproducing
kernel Hilbert spaces, e.g. splines or radial basis functions \citep[see
for instance][]{wendland05:_scatt_data_approx}, and linear prediction of
random processes, also known as Kriging \citep[see for
instance][]{chiles99}. In this paper, we shall adopt the probabilistic
framework\footnote{In fact, these two classes of methods, which have
  been studied separately, are equivalent (see for instance
  \cite{Wah70b}).}.

Thus, let us consider that an unbiased linear estimator $\hat \xi_n$ of
$\xi$ has been obtained from $\xi(x_1),\ldots,\xi(x_n)$. In particular,
we can use ordinary Kriging when the mean of $\xi(x)$ is known and
intrinsic Kriging when it is unknown, which is more often the case.

Can we expect a faster convergence when $\xi$ is replaced by $\hat
\xi_n$?  Here, we assume the computation time to evaluate $\hat
\xi_n(x)$, $x\in\XX$, conditioned on $\xi(x_i)=f(x_i)$, $i=1,\ldots,n$,
is small, which means that we can make $\abs{A_u(\hat \xi_n)}_l
-\abs{A_u(\hat \xi_n)}$ negligible with respect to $\abs{A_u(\hat
  \xi_n)} -\abs{A_u(\xi)}$. Thus, we are now interested in the
convergence of $\abs{A_u(\hat \xi_n)}$ to $\abs{A_u(\xi)}$.
Section~\ref{sec:asymptotics} shows how the convergence rate in mean
square of $\abs{A_u(\hat \xi_n)}$ to $\abs{A_u(\xi)}$ depends on the
fill distance of $\XX$ and the regularity of $\xi$. In
Section~\ref{sec:SUR}, we shall propose an algorithm to speed up this
rate by a sequential choice of the evaluation points.

\subsubsection{Intrinsic Kriging basics}
\label{sec:IK}

In this paper, we use \emph{intrinsic Kriging} (IK) to obtain a linear
predictor of $\xi$ based on a finite set of pointwise observations of
the process. We recall here the main results \citep{Mat73}. IK extends
linear prediction when the mean of $\xi(\x)$ is unknown but can be
written as a linear parametric function $m(x)=b^{ \mathsf{T}}{p}(\x)$.
Here, $p(\x)$ is a $q$-dimensional vector of base functions of a vector
space $\NN$ of translation-stable functions (in practice, all
polynomials of degree less or equal to $l$) and $b$ is a vector of
unknown parameters.  Intrinsic Kriging assumes that observed values of
$f$ are samples from a representation of an intrinsic random function
(IRF), a generalized random process defined over a space $\Lambda_l$ of
measures orthogonal to $\NN$, and characterized by its stationary
generalized covariance $k(h)$ (see the Appendix Section for more
details).

\begin{prop}[Intrinsic Kriging, \citealt{Mat73}]
Let $\xi_{\rm G}$ be an IRF$(l)$, with generalized covariance
$k({h})$. Assume $n$
observations be sample values of  the random
variables $\xi^{\mathrm{obs}}_{x_{i}}=\xi(x_{i})+N_{i}$, $i=1,\ldots
,n$, where $\xi$ is an unknown representation of $\xi_{\rm G}$ and 
the $N_{i}$s are zero-mean random variables independent of $\xi(\x)$,
with covariance matrix ${K}_N$.


The intrinsic Kriging predictor of $\xi(x)$ based on the observations,
is the linear projection $\hat \xi_n(\x) = \sum_i {\lambda }_{i,x}
\xi^{\rm obs}_{x_{i}}$ of $\xi(\x)$ onto ${\mathcal{H}}_{S}=\mathop{\rm
  span}\{\xi^{\rm obs}_{x_{i}},i=1,\ldots ,n\}$, such that the variance
of the prediction error $\xi(x)-\hat{\xi}_n(x)$ is minimized under
the constraint $\delta _{x}-\sum {\lambda }_{i,x}\delta _{x_{i}}\in
\Lambda _{l}$. The coefficients ${\lambda }_{i,x}$, $i=1,\ldots ,n$, are
solutions of a system of linear equations, which can be written in
matrix form as
\begin{equation}
\left( 
\begin{array}{cc}
{K}+{K}_{N} & {P}^{\mathsf{T}} \\ 
{P} & {0}
\end{array}
\right) \left( 
\begin{array}{c}
{{\lambda}}_{x} \\ 
{\mu}
\end{array}
\right) =\left( 
\begin{array}{c}
{k}_{x} \\ 
{p}_{x}
\end{array}
\right) \,,  \label{eq:4}
\end{equation}
where ${K}$ is the $n\times n$ matrix of generalized covariances $
k(x_{i}-x_{j})$, ${P}$ is a $q\times n$ matrix with entries ${x_j}^{i}$ for
$j=1,\ldots,n$ and multi-indexes $i=(i_1,\ldots,i_d)$  such that
$\abs{i}:=i_1+\cdots+i_d\leq l$, ${\mu}$ is a vector of Lagrange
coefficients, ${k}_{x}$ is a vector of size $n$ with entries
$k(x-x_{i})$ and ${p}_{x} $ is a vector of size $q$ with entries
$x^{i}$, $i$ such that $\abs{i}\leq l$.

The variance of the prediction error is given by 
$\sigma_n(x)^2:=\var [\xi(x)- \hat \xi_n(x)]=k(0)- \lambda_x\tr k_x -
\mu\tr {p}_x$.
\end{prop}
\begin{proof}
  See \cite{Mat73}.
\end{proof}

\subsubsection{Asymptotics}
\label{sec:asymptotics}

In this section, we shall justify that modeling the unknown $f$ by a
Gaussian random process $\xi$ and estimating $\abs{A_{u}( \xi)}$ by
$\abs{A_{u}(\hat \xi_n)}$ is well-founded.  Our objective is to
establish a mean square convergence when the evaluation points fill
$\XX$.

Classical results in approximation theory \citep[see for
instance][]{wu93:_local, light98, narcowich03:_refin,
wendland05:_scatt_data_approx} assert that the variance
$\sigma_n^2(x)$ of the IK prediction error at $x$ decreases as the
sampling density or the regularity of the covariance increases. More
precisely, if $\XX$ is a bounded domain of $\RR^d$, and the Fourier
transform of $k(h)$, $h\in \RR^d$, satisfies
$$
c_1(1+\ns{\omega}_2^2)^{-\nu} \leq \tilde k(\omega) \leq
c_2(1+\ns{\omega}_2^2)^{-\nu}\,.
$$
with $\nu>d/2$, then
\begin{equation}
  \label{eq:5}
  \ns{\sigma_n(.)}_\infty \leq C h_n^{\nu - d/2} \,,
\end{equation}
where $h_n=\sup_{y\in\XX} \min_{i} \ns{ y - x_i}_2$ is a \emph{fill distance} of
$(\x_1,\ldots,\x_n)$ in $\XX$. 

The following theorem shows that a similar result holds for the process
thresholded at a level $u$.
\begin{thm}
  Let $\xi$ be an unknown representation of an IRF($l$) $\xi_{\rm G}$, and
  $\hat \xi_n(\x)$ be   the IK predictor of $\xi$ based on
  observations $\xi(x_{i})$, $i=1,\ldots ,n$. Define $ \sigma_n(x) :=
  \var{[\xi(x) - \hat \xi_n(x)]}^{1/2}$. Then,
  $$
  \EE\Big[ (\one_{\xi(x) \geq u} - \one_{\hat \xi_n(x) \geq u})^2 \Big]
  = O(\sigma_n(x) \abs{\log (\sigma_n(x))}^{1/2}) \mbox{~~~when } \sigma_n(x)
    \rightarrow 0\,.
  $$

\end{thm}
\begin{proof}
  For all $x\in\XX$, $\xi(x) - \hat \xi_n(x)$ is Gaussian with zero-mean
  and variance $\sigma_n(x)^2$ (but is not orthogonal to $\hat \xi_n(x)$, as
  would be the case if the mean of $\xi$ were known).  Thus, $\forall
  x\in\XX$ and $\forall n\in\NNN$, we can write $\xi(x)$ as 
  \begin{equation}
    \xi(x) = (1+a_n(x)) \hat \xi_n(x) + b_n(x) + \zeta_n(x)\,,
    \label{eq:6}
  \end{equation}
  where $a_n(x), b_n(x) \in \RR$, $\zeta_n(x)$ is Gaussian and such that
  $\EE[\hat \xi_n(x) \zeta_n(x) ]=0$ and $\EE[\zeta_n(x)] = 0$.  This
  decomposition exists and is unique for every $n$.  (To simplify
  notations, from now on, we shall omit the dependence on $x$ when
  there is no ambiguity.)

  Clearly, $\var[\hat\xi_n]$ is non-decreasing and can be assumed to be
  strictly positive for $n$ large enough. Since $ \EE[ a_n \hat\xi_n]
  = -b_n $,  we have 
  \begin{equation}
    \sigma_n^2 = \var[a_n \hat\xi_n + b_n + \zeta_n] = \EE[(a_n
    \hat\xi_n + b_n + \zeta_n)^2]
    = a_n^2 \var[\hat\xi_n] + E[\zeta_n^2],
  \end{equation}
  and thus, the following upper bounds hold for $n$ large enough:
  \begin{equation}
  \left\{
    \begin{array}{l}
      \abs{a_n} \leq K_a\, \sigma_n\,,\; \abs{b_n} \leq K_b\, \sigma_n\,, \\
      \tilde \sigma_n := \EE[\zeta_n^2]^{1/2} \leq \sigma_n\,,
    \end{array}\right.
  \label{eq:7}
  \end{equation}
for some  $K_a, K_b >0$.

For some threshold $u\in\RR$, let $\alpha$ be such that 
\begin{equation}
\label{eq:8}
\alpha >  \abs{a_n u + b_n} \geq 0\,,
\end{equation}
and let $N\in \NNN$ be such that $\forall n>N$, $\abs{a_n} < 1$.
For all $n>N$, define
\begin{equation*}
\left\{
\begin{array}{lll}
h_n^- &=& \displaystyle \frac{u-b_n-\alpha}{1+a_n}\,, \\[2.5ex]
h_n^+ &=&  \displaystyle \frac{u-b_n+\alpha}{1+a_n}\,. \\
\end{array}
 \right.
\end{equation*}
Note that $h_n^- < u < h_n^+$ and that
$$
h_n^+ -h_n^-= \frac{2\alpha}{1+a_n}\,.
$$
For all $n>N$, 
 \begin{equation}
   \begin{split}
    \EE\Big[ (\one_{\xi(x) \geq u} - \one_{\hat \xi_n(x) \geq u})^2 \mid \hat
      \xi_n(x) \Big] &=
    \Psi\left(\frac{u-(1+a_n)\hat\xi_n-b_n}{\tilde\sigma_n}\right) 
    \one_{\hat \xi_n(x)<u}\\     
 &+ \Psi\left(-\frac{u-(1+a_n)\hat \xi_n-b_n}{\tilde\sigma_n}\right) \one_{\hat
      \xi_n(x)\geq u} \,,
   \end{split}
   \label{eq:9}
  \end{equation}
  in which $\Psi$ denotes the tail of the standard Gaussian distribution
  function. Since
  \begin{equation*}
      \left\{    \begin{array}{lll}
          \hat \xi_n < h_n^- & \Rightarrow &
          u-(1+a_n)\hat\xi_n-b_n > \alpha\,, \\
          \hat \xi_n > h_n^+ & \Rightarrow &
          - u +(1+a_n)\hat\xi_n+b_n > \alpha\,, 
    \end{array}\right.
  \end{equation*}
  and $\tilde \sigma_n \leq \sigma_n$, we have
  \begin{equation}
    \EE\Big[ (\one_{\xi(x) \geq u} - \one_{\hat \xi_n(x) \geq u})^2 \mid \hat
      \xi_n(x) \Big] 
      \leq \Psi\left(\frac{\alpha}{\sigma_n}\right) 
    \one_{\hat \xi_n(x) \in \RR\backslash [h_n^-, h_n^+]} + \one_{\hat
      \xi_n(x)\in [h_n^-, h_n^+]} \,.
    \label{eq:10}
  \end{equation}
  By integrating with respect to the density of $\hat \xi_n$, we
  obtain
  \begin{eqnarray}
    \EE\Big[ (\one_{\xi(x) \geq u} - \one_{\hat \xi_n(x) \geq u})^2\Big] &\leq&
    \Psi\left(\frac{\alpha}{\sigma_n}\right) + c_0 \frac{2\alpha}{1+a_n}   \\
    &\leq& \frac{\sigma_n}{\alpha \sqrt{2 \pi}}  \exp({-\frac{\alpha^2}{2\sigma_n^2}})
    +c_1 \alpha \label{eq:11}
  \end{eqnarray}
where (\ref{eq:11}) uses a standard Gaussian tail inequality. 

  The upper bound can be tighten by replacing $\alpha$ with a sequence
  $(\alpha_n)$ such that 
  $$
  \alpha_n := \sqrt{2} \sigma_n \abs{ \log (\sigma_n)}^{1/2},
  $$
  which satisfies~(\ref{eq:8}) for $n$ large enough.
  Therefore,
  \begin{equation}
    \EE\Big[ (\one_{\xi(x) \geq u} - \one_{\hat \xi_n(x) \geq u})^2\Big]
    \leq O(\sigma_n \abs{ \log (\sigma_n)}^{1/2}) \mbox{~~~when } \sigma_n
    \rightarrow 0.
  \end{equation}

\end{proof}
Hence, if $\XX$ is bounded: 
\begin{eqnarray}
  \EE\Big[(\abs{A_u(\xi)}- \abs{A_u(\hat \xi_n)} )^2\Big]
 &=& \EE\bigg[\bigg(\int_\XX \one_{\xi(x) \geq u} -
  \one_{\hat \xi_n(x) \geq u} d\mu\bigg)^2\bigg] \nonumber \\
  & \leq & \int_\XX   \EE\big[(\one_{\xi(x) \geq u} -
  \one_{\hat \xi_n(x) \geq u})^2\big] d\mu \nonumber \\
  & \leq & C  \ns{\sigma_n(.)}_\infty \abs{\log \ns{\sigma_n(.)}_\infty}^{1/2}
\end{eqnarray}
when $n\rightarrow 0$ and $\ns{\sigma_n(.)}_{\infty}\rightarrow 0$.

Therefore, this simple result shows that the mean square convergence of
$\abs{A_u(\hat \xi_n)}$ to $\abs{A_u(\xi)}$ is related to the mean
square convergence of $\hat \xi_n$ to $\xi$, hence, due to~(\ref{eq:5}),
to the regularity of the covariance and the fill-in distance of $\XX$.
Informally speaking, we can say that using an approximation will be more
efficient than a mere Monte Carlo approach if the regularity of $\xi$
compensates for the slowness of filling $\XX$, which of course increases
as the dimension $d$ of $\XX$ increases. By choosing the $x_i$s on a
lattice, the fill distance can be made such that $h_n=O(n^{-1/d})$.
Then, the convergence of $\abs{A_u(\hat \xi_n)}$ to $\abs{A_u(\xi)}$
when the $x_i$s fill $\XX$ regularly, is faster than Monte Carlo if $\nu
> 3d/2$.

\section{Convergence acceleration}
\label{sec:SUR}

\subsection{Control of convergence}

Of course, sampling $\XX$ regularly as above may be suboptimal when the
evaluations of $f$ are sequential. This section addresses the problem of
choosing a sequence $(x_n)_{n\in\NNN}$ so that the error of volume
approximation conditioned on $\xi(x_i)=f(x_i)$, $i=1,2,\ldots$ decreases
rapidly.  More precisely, a desirable strategy would consist in choosing
\begin{equation}
  \label{eq:2.1}
   x_{n} = \displaystyle \argmin_{x_{n}\in\XX} \Upsilon_n(x_{n}):= 
   \EE\big[( \abs{A_u(\xi)} - \abs{A_u(\xi_{n})})^2 \mid Z_{n-1} \big]\,,    
\end{equation}
where for all $n$, $Z_n = (\xi(x_1),\ldots,\xi(x_n))$.
Note that $\Upsilon_n(x_{n})$ can also be written as
\begin{equation}
 \Upsilon_n(x_{n}) =  \EE \big[
   \EE\big[( \abs{A_u(\xi)} -
    \abs{A_u(\xi_{n})})^2 \mid Z_{n}\big] \mid Z_{n-1} \big]\,.
 \label{eq:12}
\end{equation}
  
The distribution of $\abs{A_u(\xi)}$ conditioned on observations is
generally unknown (see \citealp{adler00excursion}, Section~4.4) and
therefore, $\EE\big[( \abs{A_u(\xi)} - \abs{A_u(\xi_{n})})^2 \mid
Z_{n}\big]$ cannot be easily determined analytically. To overcome this
difficulty, we could minimize  a Monte Carlo approximation
of~(\ref{eq:2.1}) instead, namely
\begin{equation}
  \begin{split}
   x_{n}  = &\argmin_{x_{n}\in\XX}  \Upsilon_{n,m}(x_{n}):= \\
 &\EE \bigg[  m^{-1} \sum_{i=1}^m (\abs{A_u(\xi_{n}+\zeta_{n}^i)} -
    \abs{A_u(\xi_{n})})^2  \bmid  Z_{n-1},~\{\zeta^i_{n}, {i\leq m}\}\;\bigg]\,,    
  \end{split}
  \label{eq:13}
\end{equation}
where the random processes $\zeta^i_{n}$ are $m$ independent copies of
$\xi$ conditioned on $Z_{n}=(0,\ldots,0)$. The
program~(\ref{eq:13}) becomes numerically tractable if we also replace
$\abs{A_u(\bm{\cdot})}$ by its Monte Carlo estimator
$\abs{A_u(\bm{\cdot})}_l$. Whereas simulating the conditioned processes
$\zeta^i_{n}$ is easy in principle \citep[see][chap.~7]{chiles99}, it
is also computationally intensive since it typically requires
$O(l^3)$ operations to simulate $\xi$ at given points $x_1,\ldots, x_l$.
Since $l$ has to be high enough to ensure a degree of accuracy of the estimator
$\abs{A_u(\bm{\cdot})}_l$, conditional simulations ought to be avoided.

An alternative solution is to approximate $\EE\big[( \abs{A_u(\xi)} -
\abs{A_u(\xi_{n})})^2 \mid Z_n\big]$ by $\EE\big[( \abs{A_u(\xi)}_l -
\abs{A_u(\xi_{n})}_l)^2 \mid Z_n,~\{X_i, i\leq l\}\, \big]$, for $l$ high
enough. Then,  the Minkowski inequality gives
\begin{equation}
  \label{eq:14}
  \begin{split}
    \EE\big[( \abs{A_u(\xi)}_l - &\abs{A_u(\xi_{n})}_l)^2 \mid
    Z_n,~\{X_i, i\leq l\}\, \big]^{1/2} \\
    &    \leq \frac{1}{l}\sum_{i=1}^l \EE \big[(\one_{\xi(X_i)>u} -
    \one_{\xi_n(X_i)>u})^2 \mid  Z_n,~\{X_i, i\leq l\}\,\big]^{1/2}\,.
  \end{split}
\end{equation}
This makes possible to build a \emph{stepwise uncertainty reduction}
algorithm as  presented in the next section.

\subsection{A stepwise uncertainty reduction algorithm}

Denote by $S=\{y_1,\ldots,y_l\}$ a set of $l$ independent sample values
of $X$.  Given a finite sequence $(x_i)_{1\leq i\leq n-1}$ of evaluation
points, we wish to obtain a new point $x_{n}$ that yields the largest
decrease of the upper bound of the volume approximation mean-square
error obtained in~(\ref{eq:14}), i.e.,
\begin{equation}
  x_{n}=\argmin_{x_{n} \in S} \Upsilon'_{n}(x_{n}):=\frac{1}{l}\sum_{i=1}^l \EE \big[(\one_{\xi(y_i)>u} -
  \one_{\xi_{n}(y_i)>u})^2 \mid  B_{n-1}\,\big]^{1/2} \,, 
\label{eq:15}
\end{equation}
where $B_n$ denotes the event $\{\xi(x_1) = f(x_1), \ldots,
\xi(x_n) = f(x_n)\}$, $n>0$.

A few steps are needed to transform~(\ref{eq:15}) into a numerically
tractable program. First, note that
\begin{equation}
  \begin{split}
\EE \big[(\one_{\xi(y_i)>u} -&
    \one_{\xi_{n}(y_i)>u})^2 \mid  B_{n-1}\,\big]\\ 
&  = \int_{z\in\RR} \EE\big[(\one_{\xi(y_i)>u} -
    \one_{\xi_{n}(y_i)>u})^2 \mid \xi(\x_{n}) = z, B_{n-1}\big] \\
    & \times~~p_{\xi(x_{n})\mid B_{n-1}}(z) dz\,,\quad \forall  i \in \{1,\ldots,l\},
  \end{split}
\label{eq:16}
\end{equation}
where $p_{\xi(x)\mid B_{n-1}}$ denotes the density of $\xi(x)$
conditionally to $B_{n-1}$. However, intrinsic Kriging assumes that the
mean of $\xi$ is unknown and therefore, for $x\in \XX$,
$\EE\big[(\one_{\xi(x)>u} - \one_{\xi_{n}(x)>u})^2 \mid
\xi(\x_{n})=z,B_{n-1} \big]$ cannot be determined exactly. Indeed, the
values of $a_n(x), b_n(x)$ and $\tilde \sigma_n(x)$ in~(\ref{eq:9}) are
unknown in practice.  Nevertheless,~(\ref{eq:7}) leads to the
approximation
\begin{equation}
  \label{eq:17}
  \EE\big[(\one_{\xi(x)>u} -
  \one_{\xi_{n}(x)>u})^2 \mid \xi_n(\x) \big] \approx
  \upsilon_n(x):= \Psi\left( \left\lvert \frac{u-\xi_{n}(x)}{\sigma_{n}(x)}\right\rvert\right)\,.
\end{equation}

Finally, define a discretization operator $\Delta_Q$, which
can be written for instance as
$$ \forall h\in \RR\,,\quad {\rm \Delta}_Q h = z_1
+ \sum_{i=2}^Q (z_{i} - z_{i-1}) \one_{]z_i,+\infty[}(h) $$ 
with
$z_1 < z_2 < \dots < z_Q$. We can now write~(\ref{eq:15}) as a
numerically tractable program:
\begin{multline}
x_{n} = \argmin_{x_{n}\in S} \Upsilon''_{n} (x_{n}) :=\\ \frac{1}{l}
\sum_{i=1}^l \bigg(\sum_{j=1}^Q \PR\{\Delta_Q \xi(x_{n})=z_j |\, B_{n-1}\} 
\EE\big[ \upsilon_{n}(y_i) \mid \xi(x_{n})=z_j,\;  B_{n-1}\big]\bigg)^{1/2}\,.
\label{eq:18}
\end{multline}

An informal interpretation of~(\ref{eq:18}) is that $x_{n}$ minimizes
the error of prediction of $\one_{\xi(x)>u}$ by $\one_{\xi_n(x)>u}$,
which is measured via $\upsilon_n(x)$, averaged on $\XX$ under the
distribution $\mu$, and conditioned on the observations. When
$\Upsilon''_{n} (x)$ becomes small for all $x\in S$,
$\abs{A_u(\xi_n)}_l$ conditioned on observations  provides a
good approximation of $\mathcal{P}_u$. As will be seen in
Section~\ref{sec:example}, the proposed strategy is likely to achieve very
efficient convergences.

\section{Example}
\label{sec:example}

This section provides a one-dimensional illustration of the proposed
algorithm.  We wish to estimate~(\ref{eq:1}), where $f(x)$ is a given
function defined over $\RR$ and $X\sim \mu=\mathcal{N}(0,\sigma^2)$. We
assume that $f$ is a sample path of $\xi$.  After a few iterations,
the unknown function $f$ (as shown in Figure~1) has been sampled so that
the probability of excursion $\PR\{ \xi(x) > u \mid \xi(x_i) = f(x_i),
i=1\,\ldots,n \}$ is determined accurately in the region where the
probability density of $X$ is high. This example illustrates the
effectiveness of the proposed algorithm. Note that in practice, a
parametrized  covariance has to be chosen for $\xi$ and its parameters
should be estimated from the data, using, for instance, a maximum
likelihood approach \citep[e.g.][]{Ste99}. 

\begin{figure}[htbp]
  \centering
  \includegraphics[width=0.85\textwidth]{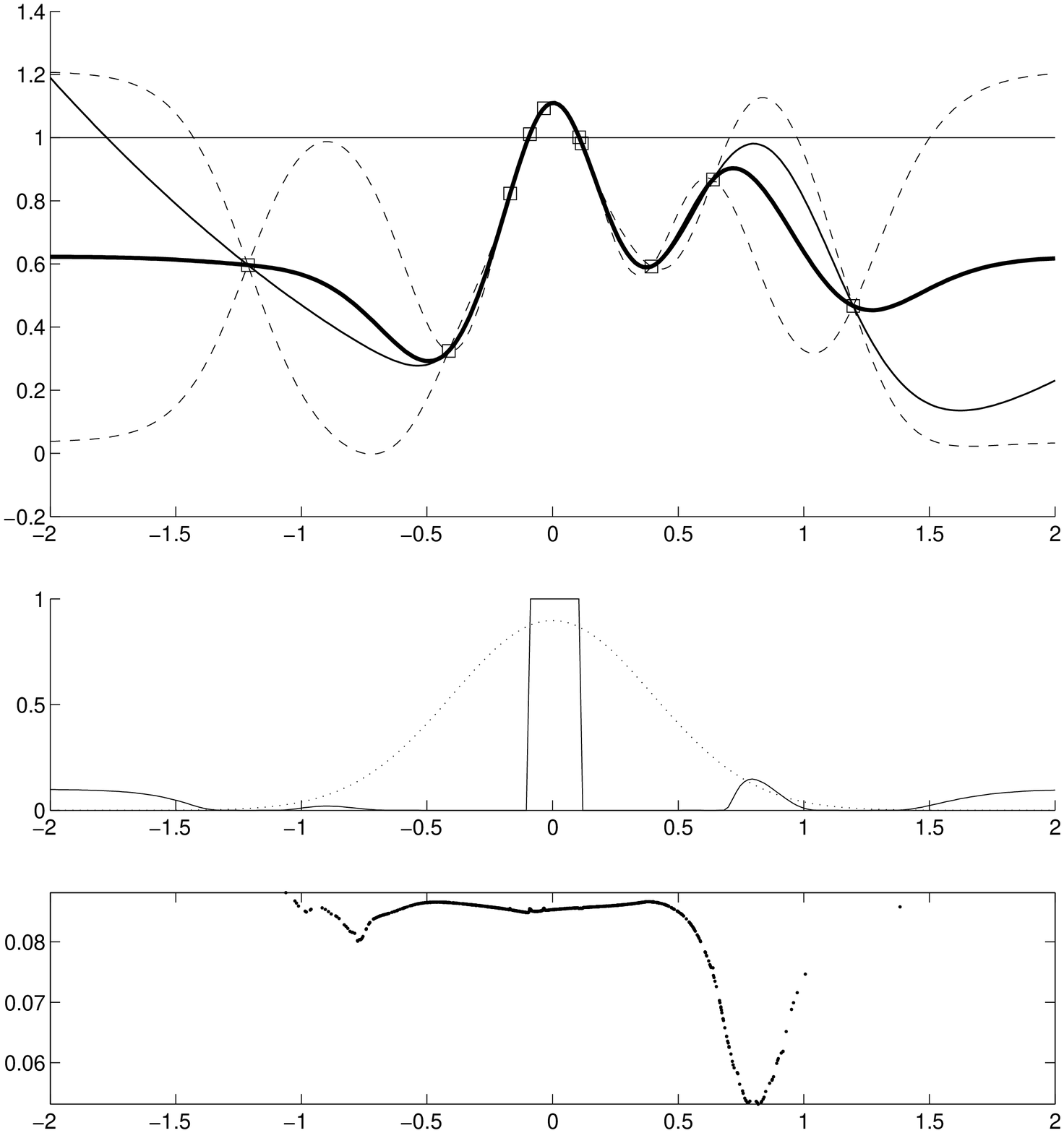}
  \captionsetup{labelsep=period}
  \caption{Top: threshold $u$ (horizontal solid line), function $f$
    (thin line), n=10 evaluations as obtained by the proposed algorithm
    using $l=800$ and $Q=20$ (squares), IK approximation $f_n$ (thick
    line), 95\% confidence intervals computed from the IK variance
    (dashed lines).  Middle: probability of excursion (solid line),
    probability density of $X$ (dotted line). Bottom: graph of $
    \Upsilon''_{n} (y_i)$, $i=1,\ldots,l=800$, the minimum of which
    indicates where the next evaluation of $f$ should be done (i.e., at
    approximately 0.75).}
  \label{fig:1}
\end{figure}

\section{Appendix : Intrinsic Random Functions}
\label{sec:IRF}

In this section, we intend to summarize the most important notions about
intrinsic random functions \citep{Mat73}.  Let ${\NN}$ be a vector space
of functions $\{{b}^{\mathsf{T}} {r}(\x),\;{b}\in {\mathbb{R}}^{l}\}$
and $\xi(\x)$ be a random process with mean $m(\x)\in {\NN}$. The main
idea of intrinsic random functions is to find some linear
transformations of $\xi(\x)$ \emph{filtering out} the mean so as to
consider a zero-mean process again.

 Let $\tilde{ \Lambda}$ be the vector space of \emph{finite-support
  measures}, \emph{i.e.}  the space of linear combinations
$\sum_{i=1}^{n}\lambda _{i}\delta _{\x _{i}}$, where $\delta _{\x}$
stands for the Dirac measure, such that for any $B\subset {\XX}$,
$\delta _{\x}(B)$ equals one if ${ \x}\in B$ and zero otherwise. Let
$\tilde{\Lambda}_{{\NN}^{\perp }}$ be the subset of the elements of
$\tilde{\Lambda}$ that vanish on~${\mathcal { N}}$. Thus, $\lambda \in
\tilde{\Lambda}_{{\NN}^{\perp }}$ implies
$$
\langle \lambda ,f\rangle :=\sum_{i=1}^{n}\lambda _{i}f(\x
_{i})=0\,,\quad \forall \,f\in {\NN}\,.
$$
In the following, we shall restrict ourselves to the case where
${\mathcal{N}}$ is a vector space of polynomials of degree at most equal
to~$l$. Denote by $\mathcal{N}_{l}$ the linear hull of all multivariate
monomials $x^i$, where $i=(i_1,\ldots,i_d)$ are multi-indexes such that
$\abs{i}:=i_1+\cdots+i_d\leq l$, and define
$\tilde{\Lambda}_{l}:=\tilde{\Lambda}_{ {\mathcal{N}}_{l}^{\perp }}$.

Let $\xi_{{\rm G}}(\lambda )$ be a linear map on $\tilde{ \Lambda}_l$,
with values in $L^{2}(\Omega,\mathcal{A},\mathsf{P})$, the space of
second-order random variables. Assume that $\EE[\xi_{{\rm G}}(\lambda
)]=0$ for all~$\lambda$ and that
$$
k(\lambda ,\mu ):= \cov[\xi_{\rm G}(\lambda),\xi_{\rm G}(\mu)] =
\sum_{i,j}\lambda _{i}\mu _{j}k(\x_{i},{ {y}}_{j})\,,
$$
where $k({ \x},{y})$ is a symmetric conditionally positive definite
function (i.e.  a function such that $k(x,y)=k(y,x)$ and
$k(\lambda,\lambda)\geq 0$ for all $\lambda \in\tilde{ \Lambda}_l$).
Then, $\xi_{ {\rm G}}(\lambda )$ is a \emph{generalized random process}
and $k(\x,{{y}})$ is called a \emph{generalized covariance} (note that
any covariance is a generalized covariance).  Let
$\tilde{\mathcal{H}}_l$ be the subspace of $L^{2}(\Omega
,\mathcal{A},\mathsf{P})$ spanned by $\xi_{{\rm G}}(\lambda )$, $\lambda
\in \tilde{\Lambda}_l$. Since random variables in
$\tilde{\mathcal{H}}_l$ are zero-mean, the inner product of
$L^{2}(\Omega ,\mathcal{A},\mathsf{P})$ can be expressed in
$\tilde{\mathcal{H}}_l$ as
$$
(\xi_{{\rm G}}(\lambda ),\xi_{{\rm G}}(\mu ))_{L^{2}(\Omega
  ,\mathcal{A}, \mathsf{P})}=k(\lambda ,\mu )\,,\quad \lambda,\mu \in
\tilde \Lambda_l\,.
$$
Thus, the bilinear form $k(\lambda ,\mu )$ endows $\tilde{ \Lambda}_l$
and $\tilde{\mathcal{H}}_{{\mathcal{N}} ^{\perp }}$ with a structure of
pre-Hilbert space. The completions $\mathcal{H}_l$ and $\Lambda _l$
of $\tilde{\mathcal{H}}_l$ and $\tilde{\Lambda}_{l}$ under this inner
product define isomorphic Hilbert spaces. $\xi_{\rm G}(\lambda)$ can
be extended on $\Lambda _l$ by continuity.
Simplifying hypotheses are introduced in the next paragraph.

Let $\tau _{{h}}:\tilde{\Lambda}_l\rightarrow \tilde{\Lambda}$ be the
translation operator such that for $ \lambda =\sum_{i}\lambda _{i}\delta
_{\x_{i}}\in \tilde{\Lambda}_l$, $\tau _{{h}}\lambda =\sum_{i}\lambda
_{i}\delta _{\x_{i}+{h}}$. Note that $\tilde{\Lambda}_l$ is stable under
translation since ${\mathcal{N}_l}$ is itself a translation-stable space
of functions. Assume further that the generalized covariance
$k(\x,{{y}})$ is invariant by translation.  In the following, we shall
write $k({h})$ with ${h} = \x - \y$ instead of $k(\x,\y)$, when the
covariance is assumed to be stationary. Then $\tau _{{h}}$ is continuous
and can be uniquely extended on $\Lambda _l $.

\begin{defn}
  Let $\xi_{{\rm G}}(\lambda )$ be a zero-mean generalized random
  process defined on $\Lambda _l$, with stationary generalized
  covariance $k({h})$. The random process ${h}\mapsto \xi_{\rm G}(\tau
  _{{h} }\lambda )$, $\lambda \in \Lambda _l$, is therefore weakly
  stationary.  $\xi_{{\rm G}}(\lambda )$, $\lambda \in \Lambda _l$, is
  then an \emph{Intrinsic Random Function} of order $l$, or IRF$(l)$ in
  short.
\end{defn}

If $\xi(\x)$, $\x\in {\XX}$, is a second-order random
process, with mean in ${\mathcal{N}_l}$ and covariance $k(x,y)$,
the linear map 
$$
\begin{array}{llcl}
\xi: & \tilde{\Lambda}_l & \rightarrow  & {\mathcal{H}}
\\ 
& \lambda =\sum_{i=1}^{n}\lambda _{i}\delta _{\x_{i}} & \mapsto  & 
\xi(\lambda ):=\sum_{i=1}^{n}\lambda _{i}\xi(\x_{i})\,,
\end{array}
$$
extends $\xi(\x)$ on $\tilde{\Lambda}_l$, where ${ \mathcal{H}}$ stands
for the Hilbert space generated by $\xi(\x)$, ${\x}\in {\XX}$.
Since $k(x,y)$ is positive definite, $(\lambda ,\mu
)_{\tilde{\Lambda}_l}:=(\xi(\lambda ),\xi(\mu ))_{{\mathcal{H}}}$
defines an inner product on $\tilde{\Lambda}_l$. Let $\Lambda _l$ be the
completion of $\tilde{\Lambda}_l$ under this inner product and extend
$\xi(\lambda )$ on $\Lambda_l$ by continuity  (a generalized random
process is thus obtained).

\begin{defn}
  Let $\xi_{{\rm G}}(\lambda )$ be an IRF$(l)$. A second-order random
  process $\xi(\x )$, $\x\in {\XX}$, is a \emph{representation} of
  $\xi_{{\rm G} }(\lambda )$ iff
$$
\xi_{{\rm G}}(\lambda )=\xi(\lambda ),\quad \forall \,\lambda \in
\Lambda_l\,.
$$

\end{defn}
If $\xi_0(\x)$ is any representation of $\xi_{\rm G}(\lambda)$, other
representations of $\xi_{\rm G}(\lambda)$ can be written as
\begin{equation}
  \label{eq:4.6}
  \xi(\x) = \xi_0(\x) + \sum_{i=1}^q B_i p_i(\x)\,,
\end{equation}
where the $p_i$s form a basis of $\NN_l$ and the $B_i$s are any
second-order random variables. Thus, the representations of an IRF$(l)$
constitute a class of random processes with mean in $\NN_l$
\citep{Mat73}.


\addcontentsline{toc}{section}{References}
\bibliography{krigeage}

\end{document}